\documentclass[12pt]{amsart}
\usepackage{amscd, amssymb,latexsym,amsmath, amscd, amsmath}
\usepackage{extarrows}
\usepackage{hyperref}
\hypersetup{%
        colorlinks,
        breaklinks=true,
        plainpages=false,%
        citecolor=blue,
        linkcolor=blue,
        urlcolor=blue,
        bookmarksopen=true,%
        bookmarksnumbered=false,%
        bookmarksdepth=5%
}
\UseRawInputEncoding
\usepackage[pagewise]{lineno}

\newcommand{\F}{\mathbb{F}}

\newcommand{\Ext}{{\rm Ext}}

\usepackage{comment}
\newcommand\restr[2]{{
  \left.\kern-\nulldelimiterspace 
  #1 
  \littletaller 
  \right|_{#2} 
  }}
\newcommand{\littletaller}{\mathchoice{\vphantom{\big|}}{}{}{}}

\def\GL{{\rm GL}}

\begin{document}

\newtheorem{theorem}{Theorem}
\newtheorem{thm}[equation]{Theorem}
\newtheorem{prop}{Proposition}
\newtheorem{cor}[equation]{Corollary}
\newtheorem{conj}{Conjecture}
\newtheorem{lemma}{Lemma}
\newtheorem{question}[equation]{Question}

\newtheorem{conjecture}[theorem]{Conjecture}

\newtheorem{example}[equation]{Example}
\numberwithin{equation}{section}

\newtheorem{remark}{Remark}

\title[On Extensions of Principal series representations]{On Extensions of Principal series representations}

\author[ Gautam H. Borisagar]{ Gautam H. Borisagar}
\author[Asfak Soneji]{Asfak Soneji}
\address{Department of Basic Sciences (Mathematics), Institute of Infrastructure Technology Research and Management(IITRAM), Maninagar (East), Ahmedabad - 380026, Gujarat, India}

\email{ghbsagar@iitram.ac.in}
\email{sonejiasfak96@gmail.com}

\begin{abstract}
  We compute $\Ext^{1}_B(\chi_1,\chi_2)$ between two characters $\chi_1,\chi_2$ of a Borel subgroup
  $B$ of a split reductive group $G$ over a finite field $\F_q,$ and make an application to the calculation of $\Ext^1_G(\pi_1,\pi_2)$ between principal series representations $\pi_1,\pi_2$
  of $G(\F_q).$ 
\end{abstract}
\maketitle

{\tiny{\textbf{Keywords:} Modular representation, extension of representations, principal series}}

\section{Introduction}\label{sec1}
Let $p$ be an odd prime and $\F_q$ be an $f$ degree extension of $\F_p.$
Let $G$ be a split reductive algebraic group defined over the finite field $\F_q$. Fix a maximal split torus $T\subseteq G$ and a Borel subgroup $B\supseteq T$ with $B = TN$, where $N$ is the unipotent radical of $B$. Note that every simple $B$-module $V$ is one-dimensional. As $N$ acts trivially on $V$, every character of $B$ can be treated as a character of $T$.
The motivation for studying extensions of representations of such groups comes from \cite{DJ} and \cite{HU}. The study of extensions of representations of certain $p$-adic groups appears in \cite{cp}, \cite{vp} and the references therein.

The main part of the article begins in Section \ref{sec2}, where we derive results about the extensions between two characters of $B$. In Proposition \ref{P3}, we prove that the Ext space between characters is nonzero if and only if characters differ by the Frobenius powers of a simple root of $T$ on $N.$

In Section \ref{sec3}, as an application, we obtain results regarding extensions between principal series representations of $G(\F_q).$ In Theorem \ref{T1}, we derive that the Ext space between principal series representations of $G(\F_q)$ is non-trivial if and only if characters differ by the Frobenius powers of simple roots.

In \cite{GJ}, a necessary and sufficient condition is provided for the extensions between certain principal series representations of $\GL_2(\F_p)$ to be split. This result can be compared with the result in Theorem \ref{T1} applied to a specific case when $G = \GL_2$ over the finite field $\F_p$.

\section{Extensions of characters of a Borel subgroup}\label{sec2}
Let $G$ be a split reductive algebraic group defined over the finite field $\mathbb{F}_{p}.$ Fix a maximal split torus $T \subseteq G$ and a Borel subgroup $B \supseteq T $ with $B=TN,$ where $N$ is the unipotent radical of $B.$ Let $\chi: T \rightarrow \mathbb{F}_{p}^\times$ be a character considered also as a character $\chi: B \rightarrow \mathbb{F}_{p}^\times.$

Recall that the action of $T$ on $N$ by conjugation, equivalently on its Lie algebra $\mathfrak{n},$ 
gives rise to the root space decomposition of $\mathfrak{n}$ as \begin{align*}
    \mathfrak{n}=\sum_{\alpha: T \rightarrow \mathbb{F}_{p}^\times}\mathfrak{n}[\alpha].
\end{align*}
In particular, there is the notion of positive roots and simple roots (of $T$ on $N$).
We have
$[N,N] \subseteq N,$ and by the Chevalley commutation relationship,   $N/[N,N]$ is a direct sum of simple root spaces $N_\alpha$:
\begin{align*}
    N/[N,N] \cong \bigoplus_{\alpha \ \text{simple}}N_{\alpha},
\end{align*}
this decomposition is as $T$ or $B$-modules.\\

We begin by proving the following result.

\begin{prop}\label{P1}
  Let $\chi_1,\chi_2$ be two characters $\chi_1,\chi_2: T \rightarrow \F_p^\times$. 
  Then, \[\Ext_{B}^{1}\left(\chi_1, \chi_2 \right)
  \neq 0 \Longleftrightarrow  \chi_1^{-1}\chi_2 = \alpha,\]
    for $\alpha$ a simple root of $T$ on $N.$ 
    Further, if $\chi_1^{-1}\chi_2 = \alpha$, for $\alpha$ a simple root of $T$ on $N,$
    treated as a character on $B,$
    \[\Ext_{B}^{1}\left(\chi_1, \chi_2 \right)
    = \mathbb{F}_{p}.\]
  
\end{prop}

\begin{proof}
  As,
  \[\Ext_{B}^{1}\left(\chi_1, \chi_2 \right) = \Ext_{B}^{1}\left({\bf{1}}, \chi_1^{-1}\chi_2 \right),\]
  where $\bf{1}$ is the trivial representation of $B$, it suffices to prove the Proposition for $\chi_1 = \bf{1}$.
  
  Thus, we suppose $\chi =  \chi_1^{-1}\chi_2: T \rightarrow \mathbb{F}_{p}^\times$ is a character, which we treat as a character of $B$.

  Suppose $\mathrm{Ext}_{B}^{1}\left(\bf{1}, \chi \right)\neq 0.$ Let us write down a nonzero element of $\mathrm{Ext}_{B}^{1}\left(\bf{1}, \chi \right)$ as an extension of $B$-modules:
    
    $${\displaystyle 0\longrightarrow \mathbb{F}_{p}[\chi] {\longrightarrow{}}E{\longrightarrow{}\mathbb{F}_{p} \longrightarrow 0}},$$
    which gives rise to a homomorphism 
    \begin{align*}
        E_{\chi} : B \longrightarrow \left ( \begin{array}{c|c}
  \chi & \ast  \\
  \hline
   & 1 \\    
\end{array} \right ) \subseteq \GL_2({\mathbb{F}_{p}}).
    \end{align*}

Restricting $E_{\chi}$ to $N \subseteq B,$ we get a homomorphism 
    \begin{align*}
     \restr{E_{\chi}}{N} : N \longrightarrow \left ( \begin{array}{c|c}
  1 & \ast  \\
  \hline
   & 1 \\    
\end{array} \right ) \subseteq \GL_2({\mathbb{F}_{p}}).
    \end{align*}
    Since the upper-triangular unipotent matrices in $\GL_2({\mathbb{F}_{p}})$ forms an abelian group,
    $\restr{E_{\chi}}{N}$ factors to give a homomorphism
\begin{align*}
     \restr{E_{\chi}}{N} : N/[N,N] \longrightarrow  \left ( \begin{array}{c|c}
  1 & \ast  \\
  \hline
   & 1 \\    
\end{array} \right ) \subseteq \GL_2({\mathbb{F}_{p}}).
    \end{align*}
    As $N/[N,N] \cong \bigoplus\limits_{\alpha \ \text{simple  roots}} \mathbb{F}_{p}[\alpha],$ we thus get a homomorphism 
      \begin{align*}
     \restr{E_{\chi}}{N/[N,N]} : \bigoplus_{\alpha \ \text{simple roots}} \mathbb{F}_{p}[\alpha]  \longrightarrow \mathbb{F}_{p}.
    \end{align*}
    The homomorphism $ \restr{E_{\chi}}{N/[N,N]}$ is equivariant under the action of $B$, giving rise to a homomorphism of $T$-modules:
    \begin{align*}
     \restr{E_{\chi}}{N/[N,N]} : \bigoplus_{\alpha \ \text{simple roots}} \mathbb{F}_{p}[\alpha]  \longrightarrow \mathbb{F}_{p}[\chi].
    \end{align*}
    This proves that the only $\chi: T \rightarrow \mathbb{F}_{p}^\times $ such that $\mathrm{Ext}_{B}^{1}\left(\bf{1}, \chi \right) \neq 0$, are the simple roots of $T$ on $N.$\\
    Conversely, given a simple root $\alpha : T \rightarrow \mathbb{F}_{p}^\times,$ and the simple root space $N_{\alpha} \subseteq N/[N,N],$ let $\psi_{\alpha} : N \rightarrow \mathbb{F}_{p}$ be a nontrivial homomorphism which is trivial on $[N,N]$ and on the simple root spaces $N_{\beta}, \beta \neq \alpha.$ In particular, $\psi_{\alpha}(t^{-1}nt)=\alpha(t^{-1})\psi_{\alpha}(n).$\\
    Define 
    $$        E_{\alpha} : B \longrightarrow  \left ( \begin{array}{c|c}
  \alpha & \ast  \\
  \hline
   & 1 \\    
\end{array}  \right )\subseteq \GL_2({\mathbb{F}_{p}})$$
to be $\psi_{\alpha} : N \rightarrow \mathbb{F}_{p}$ on $N,$ and $\alpha : T \rightarrow \mathbb{F}_{p}^\times.$ It can be seen that $E_{\alpha}$ is a homomorphism whose kernel in $B$, which we denote by $B_{\alpha},$ is a normal subgroup of $B.$
    It is easy to see that 
    \begin{align*}
       \mathrm{Ext}_{B/B_{\alpha}}^{1}\left(\bf{1}, \alpha \right) \longrightarrow \mathrm{Ext}_{B}^{1}\left(\bf{1}, \alpha \right)
    \end{align*}
    is injective. This can be seen either by considering $\Ext^{1}$ as an extension or by the inflation-restriction exact sequence in group cohomology.
    Therefore as $\mathrm{Ext}_{B/B_{\alpha}}^{1}\left(\bf{1}, \alpha \right)=\mathbb{F}_{p},\mathrm{Ext}_{B}^{1}\left(\bf{1}, \alpha \right)=\mathbb{F}_{p},$ as was to be proved.
\end{proof}
   
   The proof of the previous proposition actually proves the following stronger result.

\begin{prop}\label{P2}
  Let $B'=TN'$ be any solvable group over $\mathbb{F}_{p}$ with $T$ a torus, $N'$ the unipotent radical of $B'$ assumed to be successive extensions of ${\mathbb G}_a.$ Assume that $N'/[N',N']=\bigoplus N'[\alpha]$ where $N'[\alpha]$ are certain root spaces in $N'$ for characters $\alpha: T \rightarrow \mathbb{F}_{p}^\times.$ Then, for a character $\chi: T \rightarrow \mathbb{F}_{p}^\times, $ treated as a character of $B',$ $\mathrm{Ext}_{B'}^{1}\left(\bf{1}, \chi \right) \neq 0$ if and only if $\chi$ is one of the characters $\alpha: T \rightarrow \mathbb{F}_{p}^\times$ appearing in the decomposition
  $N'/[N',N']=\bigoplus N'[\alpha]$. Further, $\dim \Ext^1_{B'}\left(\bf{1},\alpha \right)$ is the multiplicity of the $\alpha$-eigenspace of $T$ on $N'/[N',N'].$ 
\end{prop}

\begin{remark}
    Except when $p=2,$ any homomorphism $\chi: B \rightarrow \mathbb{F}_{p}^\times$ is trivial on $N,$ thus the characters of $B$ considered here are general characters of  $B$ except when $p=2.$ 
\end{remark}
The proof of Proposition \ref{P1} can be easily modified to give a proof of the following proposition.

\begin{prop}\label{P3}
  Let $(G,B,T)$ be a triple of the split reductive group $G$ over $\mathbb{F}_{q},$ a Borel subgroup $B,$ and a maximal split torus $T$
  with $G \supseteq B \supseteq T$. Let $\chi: T \rightarrow \mathbb{F}_{q}^\times $ be a character, treated as a character of $B$.
  Then $\Ext^1_{B}\left(\bf{1}, \chi \right) \neq 0 \Longleftrightarrow \chi=\phi \circ \alpha,$ where $\alpha: T \rightarrow \mathbb{F}_{q}^\times$ is a simple root, and $\phi: \mathbb{F}_{q} \rightarrow \mathbb{F}_{q}$ is an automorphism of $\mathbb{F}_{q}$ (thus a power of the Frobenius automorphism).
\end{prop}

   The proof of Proposition \ref{P1} works except when we get a homomorphism of $T$-modules,
   \begin{align*}
     \restr{E_{\chi}}{N/[N,N]} : \bigoplus_{\alpha \ \text{simple roots}} \mathbb{F}_{p}[\alpha]  \longrightarrow \mathbb{F}_{p}[\chi].
    \end{align*}
   We will have this time a homomorphism of $T$-modules:
   \begin{align*}
     \restr{E_{\chi}}{N/[N,N]} : \bigoplus_{\alpha \ \text{simple roots}} \mathbb{F}_{q}[\alpha]  \longrightarrow \mathbb{F}_{q}[\chi],
    \end{align*}
which is not given to be $\mathbb{F}_{q}$-linear, and we need to appeal to the following lemma.

\begin{lemma}\label{L1}
  If $\chi_1$ and $\chi_2$ are any two characters of any group $A,$ thus $\chi_1: A \rightarrow \mathbb{F}_{q}^\times,$
  $\chi_2: A \rightarrow \mathbb{F}_{q}^\times$. Assume that $\chi_1: A \rightarrow \mathbb{F}_{q}^\times$ is surjective.
  Let $\mathbb{F}_{q}[\chi_1]$, resp. $\mathbb{F}_{q}[\chi_2]$ denote the space $\mathbb{F}_{q}$ on which $A$ operates by $\chi_1$, resp. $\chi_2$.  
Then $\mathbb{F}_{q}[\chi_1] \cong \mathbb{F}_{q}[\chi_2] $ as $\mathbb{F}_{p}[A]$-modules if and only if $\chi_1=\phi(\chi_2)$ where $\phi: \mathbb{F}_{q} \rightarrow \mathbb{F}_{q}$ is an isomorphism of fields, thus a power of the Frobenius.
\end{lemma}
\begin{proof}
    As $\mathbb{F}_{q}[\chi_1]$ (resp., $\mathbb{F}_{q}[\chi_2]$) denotes $\mathbb{F}_{q}$ with an action of $A$ via $\chi_1$ (resp., $\chi_2$), suppose $\mu: \mathbb{F}_{q}[\chi_1] \xlongrightarrow{\sim} \mathbb{F}_{q}[\chi_2],$ as $A$-modules. Then $\mu$ is an $\mathbb{F}_{p}$-linear isomorphism of $\mathbb{F}_{q}$ with $\mathbb{F}_{q}$ such that:
    \begin{equation}\label{E1}
        \mu(\chi_1(a)f)=\chi_2(a)\mu(f), \  \forall a \in A, f\in \mathbb{F}_{q}. 
    \end{equation}
    Assume by scaling that $\mu(1)=1.$ Then by equation \ref{E1}, we have $$\mu(\chi_1(a))=\chi_2(a), \ \ \forall a \in A. $$
    Therefore for all $a,b \in A$, we have:
    \[\mu(\chi_1(a)\chi_1(b))=\mu(\chi_1(ab))=\chi_2(ab)=\mu(\chi_1(a))\mu(\chi_1(b)).\]

    As we have assumed that $\chi_1: A \rightarrow \mathbb{F}_{q}^\times$ is surjective, we find that $\mu: \mathbb{F}_{q} \rightarrow \mathbb{F}_{q}$ preserves both additive and multiplicative structures, and hence is a field automorphism of $\mathbb{F}_{q}.$ Thus $\mu$ is a power of the Frobenius automorphism $x\rightarrow x^p$.
\end{proof}

\section{On Extensions of Principal series representations}\label{sec3}
The computation of $\mathrm{Ext}_{B}^{1}\left(\bf{1}, \chi \right)$ done in the previous section
has, as a consequence, the following theorem.

\begin{theorem}\label{T1} 
Let $(G,B,T)$ be a triple consisting of a split reductive group $G$ over ${\mathbb F}_q$, a Borel subgroup $B$ containing a maximal torus $T$.
    Let $\chi_1,$ $\chi_2 : T \rightarrow \mathbb{F}_{q}^\times$ be two characters on a maximal split torus $T,$ treated as characters on $B.$ Let  $\Pi_{\chi_{1}}= \operatorname{Ind}_{B}^{G} (\chi_1)$ and $\Pi_{\chi_{2}}= \operatorname{Ind}_{B}^{G} (\chi_2)$ be the corresponding principal series representations of $G=G(\mathbb{F}_{q}).$ Then, 
    \[\mathrm{Ext}_{G}^{1}\left(\Pi_{\chi_{1}}, \Pi_{\chi_{2}} \right) \neq 0 \implies \chi_{1}^{-1}\chi_{2}^w=\phi \circ \alpha, \]
for some $w\in W,$ an element in the Weyl group, and where $\alpha : T \rightarrow \mathbb{F}_{q}^\times$ is a simple root, and $\phi: \mathbb{F}_{q} \rightarrow \mathbb{F}_{q}$ is an automorphism of $\mathbb{F}_{q}$
    (thus a power of the Frobenius automorphism).
Further, if $\alpha$ is a simple root and $\chi_{1}^{-1}\chi_{2}=\phi \circ \alpha$ for $\phi$ an automorphism of $\F_{q},$ then $\mathrm{Ext}_{G}^{1}\left(\Pi_{\chi_{1}}, \Pi_{\chi_{2}} \right) \neq 0.$
    
\end{theorem}
\begin{proof}
  The proof of this theorem depends on the $\mathrm{Ext}$ version of Frobenious reciprocity according to which,
  \[\mathrm{Ext}_{G}^{i}\left(\operatorname{Ind}_{B}^{G} \chi_1, \operatorname{Ind}_{B}^{G} \chi_2 \right) \cong \mathrm{Ext}_{B}^{i}\left( \chi_1, \restr{\operatorname{Ind}_{B}^{G} \chi_2}{B} \right) {\rm~~ for ~~all~~} i \geq 0.\]

    By the Bruhat decomposition, $\operatorname{Ind}_{B}^{G} \chi_2$, as a $B$-module, is a direct sum:
    $$\restr{\operatorname{Ind}_{B}^{G} (\chi_2)}{B} \cong \bigoplus_{w \in W} \operatorname{Ind}_{B \cap B^{w}}^{B} (\chi_2^{w}).$$
    Therefore, $$\mathrm{Ext}_{B}^{i}\left( \chi_1, \restr{\operatorname{Ind}_{B}^{G} \chi_2}{B} \right)\cong \mathrm{Ext}_{B}^{i}\left( \chi_1,  \bigoplus_{w \in W} \operatorname{Ind}_{B \cap B^{w}}^{B} (\chi_2^{w}) \right),$$
    $$\hspace{2cm} \cong \bigoplus_{w \in W} \mathrm{Ext}_{B\cap B^{w}}^{i}\left(\chi_1, \chi_2^{w}  \right).$$
    The subgroup $B_w'=B\cap B^{w} \subseteq B$ contains $T,$ and is of the form for which Proposition \ref{P3} holds. Therefore,
    $\mathrm{Ext}_{B\cap B^{w}}^{1}\left(\chi_1, \chi_2^{w}  \right)$ is nonzero if and only if $\chi_{1}^{-1}\chi_{2}^{w}$ is one of the roots of $T$ in $N'/[N',N']$ for $N'$ the unipotent radical of $B_w'=B\cap B^{w}.$ It is not clear what are the characters of $T$ appearing in $N'/[N',N']$ for a general choice of $w \in W,$ except that for $w=1,$ these are the simple roots on $T$ on $B,$ and that for a general $w$, these are {\em certain} roots of $T$ on $B,$ proving the theorem.
    \end{proof}
Frobenius reciprocity, and the results of the previous section, also prove the following result much more directly.

\begin{prop} Let $(G,B,T)$ be a triple consisting of a split reductive group $G$ over ${\mathbb F}_q$, a Borel subgroup $B$ containing a maximal torus $T$.
  Let $\chi_1: G \rightarrow \mathbb{F}_{q}^\times$ be a character on $G$, and $\chi_1|_T$ its restriction to $T$.
  Further, let $\chi_2 : T \rightarrow \mathbb{F}_{q}^\times$ 
  be a character of $T,$ treated as a character of $B$,
  and $\Pi_{\chi_{2}}= \operatorname{Ind}_{B}^{G} (\chi_2)$ be the corresponding principal series representations of $G=G(\mathbb{F}_{q}).$ Then, 
  \[ \mathrm{Ext}_{G}^{1}\left(\chi_1, \Pi_{\chi_{2}} \right) \neq 0 \Longleftrightarrow  (\chi_{1}^{-1}\chi_{2}^w)|_T =\phi \circ \alpha,\]
  for some $w\in W,$ an element in the Weyl group, and where $\alpha : T \rightarrow \mathbb{F}_{q}^\times$ is a simple root, and $\phi: \mathbb{F}_{q} \rightarrow \mathbb{F}_{q}$ is an automorphism of $\mathbb{F}_{q}$
    (thus a power of the Frobenius automorphism).
  Further, if $\alpha$ is a simple root and $(\chi_{1}^{-1}\chi_{2})|_T=\phi \circ \alpha$ for $\phi$ an
  automorphism of $\mathbb{F}_{q},$ then $\mathrm{Ext}_{G}^{1}\left({\chi_{1}}, \Pi_{\chi_{2}} \right) = \F_q.$
    
\end{prop}

\noindent{\bf  Acknowledgements.}
The authors express their gratitude to Dipendra Prasad for his
valuable guidance and help with this work. The authors thank Eknath Ghate for his suggestions.
Additionally, the second author acknowledges the financial assistance as a Junior Research Fellowship (NET) from UGC, India.

\vspace{1cm}

\bibliographystyle{alpha}

\end{document}